\documentclass[12pt,draftcls,onecolumn]{IEEEtran}

\usepackage{amsfonts,amsmath,amscd,graphicx,hyperref}

\newtheorem{lem}{Lemma}
\newtheorem{defn}{Definition}

\newcommand{\YY}{{\mathcal Y }}

\newcommand{\UU}{{\mathcal  U}}
\newcommand{\LL}{{\mathcal  L}}
\newcommand{\GG}{{\mathfrak  g}}

\newcommand{\RR}{{\mathbb R}}

\newcommand{\dv}[2]{{\frac{\partial #1}{\partial #2}}}

\newcommand{\dotex}{{\frac{d}{dt}}}

\newcommand{\vw}{{\vec{\omega}}}

\newcommand{\va}{\vec{a}}

\newcommand{\vG}{\vec{G}}
\newcommand{\vB}{\vec{B}}

\usepackage{xcolor}

\begin{document}

\title{Non-linear Symmetry-preserving Observers on Lie Groups}

\author{Silvere Bonnabel, Philippe Martin and  Pierre Rouchon
\thanks{
Centre Automatique et Syst\`{e}mes,  Ecole des Mines de Paris, 60
boulevard Saint-Michel,75272 Paris CEDEX 06, FRANCE {\tt\small bonnabel@montefiore.ulg.ac.be}, {\tt\small philippe.martin@ensmp.fr},
{\tt\small pierre.rouchon@ensmp.fr}}
}
\maketitle

\begin{abstract}
In this paper we give a geometrical framework for the design of
observers on finite-dimensional Lie groups for systems which possess
 some specific symmetries. The design and the error (between true and estimated state)
 equation are explicit and intrinsic.  We consider also a particular case: left-invariant
systems on Lie groups with right equivariant
 output. The theory yields a class of observers such that error
equation is autonomous. The observers converge locally around
\emph{any} trajectory,
 and the global behavior is independent from the trajectory, which reminds the linear
  stationary case.
\end{abstract}
%\begin{keywords}
%Lie group, symmetries, invariance, nonlinear asymptotic observers,
%inertial navigation.
%\end{keywords}

\section{Introduction}

Symmetries (invariances) have been used to design controllers and for
optimal control theory
(\cite{fagnani-willems-siam93,koon-marsden-siam97,grizzle-marcus-ieee85,respondek-tall-scl02,martin-et-al-cocv03,morin-ieee03}),
 but
 far less for the design of observers. \cite{arxiv-07} develops a theory
of symmetry-preserving observers and presents three non-linear  observers for
 three examples of engineering interest:  a chemical reactor, a
 non-holonomic car, and an inertial navigation system. In the two latter examples the state
 space and the group of symmetry have the same dimension and (since the action is free)
the state
 space can be identified with the group (up to some discrete group). Applying the general theory
 to the Lie group case, we develop here a proper theory of symmetry-preserving observers
 on Lie groups.
The advantage over \cite{arxiv-07}
 is that the observer design is explicit (the implicit function theorem is not needed)
and intrinsic, the error equation and its first-order approximation
can be computed explicitly, and are intrinsic, and all the formulas
are globally defined. Moreover, this paper is a step further in the
symmetry-preserving observers theory since \cite{arxiv-07} does not
deal at all with convergence issues in the general case. Here using
the explicit error equation we introduce a new class of trajectories
 around which we build convergent observers. In the case of section
 \ref{lie2:sec} a class of first-order convergent observers around any trajectory is
given. The theory applies to various systems of engineering interest
modeled as invariant systems on Lie groups, such as
 cart-like vehicles and rigid bodies in space.
In particular it is well suited to attitude
 estimation and some inertial navigation examples.

The paper is organized as follows: in section \ref{lie:sec} we give a general
 framework for symmetry-preserving observers on Lie
groups. It explains the general form
 of the observers \cite{mahony-et-al-IEEE,hamel-mahony-icra06,bonnabel-rouchon-LN05} and
\cite{arxiv-07}
based on
 the group structure of $\mbox{SO}(3)$ and (resp.) $\mbox{SE}(2)$, without considering the
convergence issues. The design, the error equation and its first-order approximation are given explicitly. It is theoretically explained why the error
equation in the car example of \cite{arxiv-07} does not depend on the
trajectory (although it
 depends on the inputs). Then we introduce a new class  of trajectories called
 permanent trajectories which extend the notion of equilibrium point for
systems with symmetries: making
a symmetry-preserving observer around such a trajectory boils down to
make a Luenberger observer around an equilibrium point. We characterize permanent trajectories
geometrically and give a
 locally convergent observer around any permanent trajectory.

In section \ref{lie2:sec} we consider the special case of a left-invariant system
 with right equivariant output.
It can be looked at as the motion of a generalized rigid body in
space with measurements expressed in the body-fixed frame, as it will be explained
 in section \ref{right-left:sec}. Thus it applies to some inertial navigation examples. In particular
 it allows to explain theoretically why the error equation in the inertial navigation
 example of \cite{arxiv-07} is autonomous.
 A class of first-order
convergent observers such that the error equation is autonomous is
derived. This property reminds much of the linear stationary case.
We also explore the links between right equivariance of the output
map and observability.

Preliminary versions of section \ref{lie2:sec} can be found in
\cite{bonnabel:cifa06,Bonnabel-et-al-lie-ifac08}.

 \section{Symmetry-preserving observers on Lie groups}\label{lie:sec}

\subsection{Invariant observer and error equation}
Consider the following system :
\begin{align}
\frac {d}{dt} x(t) &= f(x,u)\label{dyng1:eq}\\
y&=h(x,u)\label{output:eq}\end{align} where $x\in G$,
$u\in\UU=\RR^m$, $y\in\YY= \RR^p$ (the whole theory can be easily
adapted to the  case where $\UU$ and $\YY$ are smooth $m$ and
$p$-dimensional manifolds, for instance Lie groups), and $f$ is a
smooth vector field on $G$. $u\in\UU$ is a known input (control,
 measured perturbation, constant parameter, time $t$ etc.).
\begin{defn}\label{group:action:def}
Let $G$ be a Lie Group with identity~$e$ and $\Sigma$ an open set
(or more generally a manifold). A \emph{left} group action
$(\phi_g)_{g\in G}$ on~$\Sigma$ is a smooth map
\[ (g,\xi)\in G\times\Sigma\mapsto\phi_g(\xi)\in\Sigma \]
such that:
\begin{itemize}
 \item $\phi_e(\xi)=\xi$ for all~$\xi$
 \item $\phi_{g_2}\bigl(\phi_{g_1}(\xi)\bigr)=\phi_{g_2g_1}(\xi)$
 for all $g_1,g_2,\xi$.
\end{itemize}
\end{defn}
In analogy one defines a right group action the same way except that $\phi_{g_2}\bigl(\phi_{g_1}(\xi)\bigr)=\phi_{g_1g_2}(\xi)$
 for all $g_1,g_2,\xi$.
Suppose $G$ acts on the left on
 $\UU$ and $\YY$  via
   $\psi_g : \UU\rightarrow \UU$ and $\rho_g : \YY\rightarrow \YY$.
 Suppose the
dynamics \eqref{dyng1:eq} is invariant in the sense
of \cite{arxiv-07} where the group action on the state space (the group
itself) is made of left multiplication: for any $g\in G$, $
 DL_g f(x,u)= f(gx,\psi_g(u))
$, i.e: $$ \forall ~x,g \in G~~ f(L_g(x),\psi_g(u))=DL_{g}f(x,u)$$
where $L_g : x\mapsto gx$ is the left multiplication on $G$, and
$DL_{g}$ the induced map on the tangent space. $DL_{g}$ maps the
tangent space $TG|_x$ to $TG|_{gx}$. Let  $R_{g}: x\mapsto xg$
denote the right multiplication and $DR_{g}$ its induced map on the
tangent space. As in \cite{arxiv-07}, we suppose that the output
$y=h(x,u)$  is equivariant, i.e,  $
h\bigl(\varphi_g(x),\psi_g(u)\bigr)=\rho_g\bigl(h(x,u)\bigr)$ for
all $g,x,u$.
\begin{defn}\label{inv:def}
Consider the change of variables $X=gx$, $U=\psi_g(u)$ and
$Y=\rho_g(y)$. The system \eqref{dyng1:eq}-\eqref{output:eq} is
left-invariant with equivariant output if for all $g\in G$ it is
unaffected by the latter transformation: $\frac {d}{dt} X(t) =
f(X,U),~ Y=h(X,U)$.
\end{defn}
We are going to build observers which respect the symmetries
(left-invariance under the group action) adapting the constructive
method of \cite{arxiv-07} to the Lie group case.
\subsubsection{Invariant pre-observers}
Following \cite{olver-book99} (or \cite{arxiv-07}) consider the action
$(\phi_g)_{g\in G}$
of $G$ on $\Sigma=\RR^s$ where $s$ is any positive integer.
Let $(x,z)\in G\times\RR^s$, one can compute (at most) $s$
functionally independent scalar
invariants of the variables $(x,z)$ the following way: $I(x,z)=\phi_{x^{-1}}(z)\in\RR^s$ .
It has the property that any invariant real-valued function $J(
x,z)$ which verifies $J(gx,\phi_g(z))=J(x,z)$ for all $g,x,z$ is a
function of the components of $I(x,z)$:  $J(x,z)=\mathcal{H}(I(x,z))$.
Applying this general method we find
 a complete set of invariants of $(x,u)\in G\times\UU$:
\begin{equation}\label{LieI:eq}
I(x, u)= \psi_{x ^{-1}} (u) \in \UU.
\end{equation}

Take $n$ linearly
independent vectors $(W_1, \ldots, W_n)$ in $TG|_e=\GG$, the Lie
algebra of the group $G$. Define $n$ vector fields by the invariance relation $w_i(x)=DL_x W_i\in TG|_x,
i=1...n,x\in G$. The vector fields form an invariant frame \cite{olver-book99}. According to \cite{arxiv-07}
\begin{defn}[pre-observer]\label{preobs:def}
The  system $
 \dotex{\hat{x}} = F(\hat{x},u,y)
$ is a \emph{pre-observer} of~\eqref{dyng1:eq}-\eqref{output:eq} if
$F\bigl(x,u,h(x,u)\bigr)=f(x,u)$ for all $(x,u)\in G\times \UU$.
\end{defn}
The definition does not deal with convergence; if moreover
$x(t)^{-1}\hat x(t)\rightarrow e$ as $t\rightarrow+\infty$ for every
(close) initial conditions, the pre-observer is an (asymptotic)
\emph{observer}. It is said to be  is \emph{G-invariant} if $
F\bigl(g\hat x,\psi_g(u),\varrho_g(y)\bigr)
 = DL_g F(\hat x,u,y)$ for all
$(g,\hat x,u,y)\in G\times G \times \UU\times\YY$.
\begin{lem}Any invariant pre-observer reads
\begin{equation}\label{LieIO:eq}
 \dotex \hat x = f(\hat x , u) +
 DL_{\hat x}   \left(
\sum_{i=1}^{n} \LL_i\left(\psi_{\hat x ^{-1}} (u),
    \rho_{\hat x ^{-1}}(y)\right)W_i
 \right)
\end{equation}
\end{lem}
where the $\LL_i$ are any smooth functions of their arguments such
that $\LL_i\left(\psi_{\hat x ^{-1}} (u),
    h(e,\psi_{\hat x ^{-1}} (u))\right)=0$. The proof is analogous to \cite{arxiv-07}: one can write
$DL_{\hat x^{-1}}(\dotex \hat x-f(\hat x,u))=\sum_{i=1}^{n}
 \mathcal{F}_i\left(\hat x,u,y\right)W_i\in \GG$, where the $\mathcal{F}_i's$ are
 invariant scalar functions of their arguments. But a complete set of invariants of
 $\hat x,u,y$ is made of the components of $(\psi_{\hat x ^{-1}} (u),
    \rho_{\hat x ^{-1}}(y))$, thus $\mathcal{F}_i\left(\hat x,u,y\right)=\LL_i\left(\psi_{\hat x ^{-1}} (u),
    \rho_{\hat x ^{-1}}(y)\right)$. And when $\hat x=x$ we have
    $\rho_{\hat x ^{-1}}(y)=h(\hat x^{-1} x , \psi_{\hat x ^{-1}}
     (u))= h(e,\psi_{\hat x ^{-1}} (u))$ and the $\LL_i$'s cancel.

\subsubsection{Invariant state-error dynamics}
Consider the invariant state-error $\eta=x^{-1} \hat x\in~G$. It is invariant by left multiplication
: $\eta=(gx)^{-1} (g\hat x)$ for any $g\in G$. Notice
that a small error corresponds to $\eta$ close to $e$. Contrarily to
\cite{arxiv-07}, the time derivative of $\eta$ can be computed
explicitly. We recall $R_g$ denotes the right multiplication map on
$G$. Since we have
\begin{itemize}
  \item for any $g_1,g_2\in G$, $DL_{g_1} DL_{g_2}=DL_{g_1g_2}$,
$DR_{g_1} DR_{g_2}=DR_{g_2g_1}$, $DL_{g_1} DR_{g_2}=DR_{g_2}
DL_{g_1}$
  \item $I(\hat x, u) = \psi_{\hat x^{-1}} (u)=\psi_{(x\eta)^{-1}} (u)$
  \item $\rho_{\hat
x^{-1}} (h(x,u))= h(\hat x^{-1} x , \psi_{\hat x ^{-1}} (u))$ writes
$ \rho_{\hat x^{-1}} (y)= h(\eta^{-1}, \psi_{(x\eta)^{-1}}(u)) $
\item $\dotex \eta=\dotex ( x^{-1} \hat x) = DL_{x^{-1}} \dotex \hat x - DR_{\hat x}
\dotex
  x^{-1}$ with $\dotex x^{-1}= -DL_{x^{-1}} DR_{x^{-1}} \dotex
  x$
\end{itemize}
the error dynamics reads
\begin{multline}\label{er:lie:eq}
   \dotex \eta =  DL_\eta   f(e,\psi_{(x\eta)^{-1}} (u))
   - DR_\eta   f(e,\psi_{x^{-1}}(u))
  \\ + DL_\eta
   \left(
\sum_{i=1}^{n} \LL_i\left(\psi_{(x\eta)^{-1}} (u),
    h(\eta^{-1},\psi_{(x\eta)^{-1}} (u)) \right)W_i
 \right)
 .
\end{multline}
The invariant error $\eta$ obeys a differential equation that is
coupled to the system trajectory $t\mapsto(x(t),u(t))$  only via the
invariant term $I(x,u)=\psi_{x^{-1}} (u)$. Note that when
$\psi_g(u)\equiv u$ the invariant error dynamics is independent of
the state trajectory $x(t)$! This the reason why we have this
property
 in the non-holonomic car example of \cite{arxiv-07}.

\subsubsection{Invariant first order approximation}
For $\eta$ close to $e$, one can set in~\eqref{er:lie:eq}
$\eta=\exp(\epsilon\xi)$ where $\xi$ is an element of the Lie
algebra $\GG$ and $\epsilon\in\RR$ is small. The linearized
invariant state error equation can always be written in the same
tangent space $\GG$: up to order second terms in $\epsilon$
\begin{multline}
\label{er:lie:lin:eq}
 \dotex \xi = [\xi , f(e,\psi_{x^{-1}}(u))]
 - \dv{f}{u}(e,\psi_{x^{-1}}(u)) \dv{\psi}{g}(e,\psi_{x^{-1}}(u))
 \xi
 \\
 -\sum_{i=1}^{n} \left(\dv{\LL_i}{h}(\psi_{x^{-1}}(u),h(e,\psi_{ x ^{-1}} (u)) \dv{h}{x}(e,\psi_{x^{-1}}(u)) \xi\right) W_i
\end{multline}
where [,] denotes the Lie bracket of $\GG$, $\psi$ is viewed as a function of $(g,u)$,
and $\dv{\LL_i}{h}$ denotes the
partial derivative of $\LL_i$ with respect to its second argument.  The gains
$\dv{\LL_i}{h}(\psi_{x^{-1}}(u),h(e,\psi_{ x ^{-1}} (u))$ can be tuned via linear techniques to achieve
local convergence.

\subsection{Local convergence around permanent trajectories}\label{permtraj:ssec}
The aim of this paragraph is to extend local convergence results
around an equilibrium point  to a class of trajectories
we call permanent trajectories.
\begin{defn}\label{permtraj:def}
A trajectory of~\eqref{dyng1:eq} is permanent if~
$I(x(t),u(t))=\bar I$ is independent of $t$.
\end{defn}Note that adapting this definition to the general case of symmetry-preserving
 observers \cite{arxiv-07} is straightforward.
Any trajectory of the system verifies
$
\frac{d}{dt}x(t)=DL_{x(t)}f(e,\psi_{x(t)^{-1}}(u(t)))
$
 thanks to the
invariance of the dynamics. It is permanent if
$I(x(t),u(t))=\psi_{x^{-1}(t)}(u(t))=\bar u$ is independent of $t$.
The permanent trajectory $x(t)$ is then given by $x(0)\exp(t \bar w)
$ where $\bar w$ is the left invariant vector field associated to
$f(e,\bar u)$. Thus $x(t)$ corresponds, up to a left translation
defined by the initial condition,  to a one-parameter sub-group.

Let us make an observer around an arbitrary permanent trajectory:
denote by $(x_r(t),u_r(t))$ a permanent trajectory associated to
$\bar u =\psi_{\bar x_r^{-1}(t)}{u_r(t)}$. Let us suppose we made an
invariant observer following~\eqref{LieIO:eq}. Then the error
equation~\eqref{er:lie:eq} writes
\begin{multline}\label{er:perm:lie:eq}
   \dotex \eta =  DL_\eta   f(e,\psi_{\eta^{-1}} (\bar u)) - DR_\eta   f(e,\bar u)
  + DL_\eta
   \left(
\sum_{i=1}^{n} \LL_i\left(\psi_{\eta^{-1}} (\bar u ),
   h(\eta^{-1},\psi_{\eta^{-1}} (\bar u )) \right)W_i
 \right)
 .
\end{multline}
since $\psi_{(x_r\eta)^{-1}}(u)=\psi_{\eta^{-1}}(\psi_{
x_r^{-1}}(u))=\psi_{\eta^{-1}}(\bar u)$. The first order
approximation~\eqref{er:lie:lin:eq} is now a time invariant system:
\begin{equation*}
 \dotex \xi = [\xi , f(e,\bar u)]
 - \dv{f}{u}(e,\bar u) \dv{\psi}{g}(e,\bar u)
 \xi
-\sum_{i=1}^{n} \left(\dv{\LL_i}{h}(\bar u,h(e,\bar u)) \dv{h}{x}(e,\bar u) \xi\right) W_i
\end{equation*}
Let us write $\xi$ and $f(e,u)$ in the frame defined by the $W_i$'s:
$\xi=\sum_{k=1}^{n} \xi^k W_k$ and $f(e,\bar u)=\sum_{k=1}^{n} \bar
f^k W_k$. Denote by $C_{ij}^k$ the structure constants associated
with the Lie algebra of $G$: $[W_i,W_j] =\sum_{k=1}^{n} C_{ij}^k
W_k$. The above system reads:
\begin{equation}\label{LinPermObs:eq}
 \dotex \xi
 = (A  + \bar\LL C)
 \xi
\end{equation}
where
\begin{multline*}
   A = \left(  \sum_{k=1}^{n} C_{jk}^i \bar f^k-\left[\dv{f}{u}(e,\bar u) \dv{\psi}{g}(e,\bar u)\right]_{i,j}
   \right)_{1\leq i,j\leq n}, \\
   \bar\LL = \left(- \dv{\LL_i}{h_k}(\bar u, h(e,\bar u))\right)_{\tiny
       \begin{array}{c} 1\leq i\leq n \\ 1\leq k\leq p\end{array}
       }
       ,\quad
   C = \left( \dv{h_k}{x_j}(e, \bar u)\right)_{\tiny
       \begin{array}{c} 1\leq k\leq p \\ 1\leq j\leq n\end{array}
       }
\end{multline*}
where $(x_1,\ldots,x_n)$ are the local coordinates   around $e$
defined by the exponential map:  $x=\exp(\sum_{i=1}^n x_i W_i)$. If
we assume that the pair $(A,C)$ is observable  we can choose the
poles of $A+\bar \LL C$ to get an invariant and locally convergent
observer around any permanent trajectory associated to
 $\bar u$. Let $W(x)=[W_1(x),..,W_n(x)]$. It suffices to take:
\begin{equation}\label{LiePermIO:eq}
 \dotex \hat x = f(\hat x , u(t)) +
W(\hat x)\bar\LL \rho_{\hat x^{-1}}(y(t))
\end{equation}

\subsubsection*{Examples}
In the non-holonomic car example of \cite{arxiv-07}, permanent
trajectories are made of lines and circle with constant speed. In
the inertial navigation example of \cite{arxiv-07},
$\psi_{x^{-1}}(u)=\begin{pmatrix}
  q\ast\omega\ast q^{-1}\\ q\ast(a+v\times\omega)\ast q^{-1}
\end{pmatrix}$, a trajectory is permanent if $q\ast\omega\ast
q^{-1}$ and $q\ast(a+v\times\omega)\ast q^{-1}$ are independent of
$t$. Some computations show that any permanent trajectory reads:
\begin{align*}
    q(t)&= \exp\left(\frac{\Omega}{2} t\right) \ast q_0
    \\
    v(t)&= q_0^{-1}\ast \left((\lambda\Omega t + \Upsilon +
    \exp\left(-\frac{\Omega}{2} t\right)\ast \Gamma\ast \exp\left(\frac{\Omega}{2} t\right)
     \right) \ast q_0
\end{align*}
where $\Omega$, $\Upsilon$ and $\Gamma$ are constant vectors of
$\RR^3$, $\lambda$ is a constant scalar and $q_0$ is a unit-norm
quaternion. Theses constants can be arbitrarily chosen. Hence, the
general permanent trajectory corresponds, up to a Galilean
transformation, to an helicoidal motion uniformly accelerated along
the rotation axis when $\lambda \neq 0$; when $\lambda$ tends to
infinity and $\Omega$ to $0$, we recover as a degenerate case a
uniformly accelerated line. When $\lambda=0$ and $\Gamma=0$ we
recover a coordinated turn.

\section{Left invariant dynamics and right equivariant
output}\label{lie2:sec}

\subsection{Invariant observer and error equation}
\subsubsection{Left invariant dynamics and right equivariant output}\label{right-left:sec}
Consider the following system:
\begin{align}
\frac {d}{dt} x(t) &= f(x,t)\label{dyng}\\
y&=h(x)\label{output1:eq}\end{align} where we still have $x\in G$,
$y\in\YY$, and $f$ is a smooth vector field on $G$. Let us suppose
the dynamics \eqref{dyng} is \emph{left}-invariant (see e.g
\cite{arnold-book-3}), i.e:
 $ \forall g,x \in G~~ f(L_g(x),t)=DL_{g}f(x,t)$. For all  $g\in G$, the
transformation $X(t)=g  x(t)$ leaves the dynamics equations
unchanged: $\frac{d}{dt}X(t)=f(X(t),t)$. As in \cite{arnold-book-3}
let $\omega_s=DL_{x^{-1}}\dotex x \in \mathfrak g$. Indeed one can
look at any left invariant dynamics on $G$ as a motion of a
``generalized rigid body" with configuration space $G$. Thus one can
look at $\omega_s(t)=f(e,t)$ as the ``angular velocity in the body",
where $e$ is the group identity element (whereas $DR_{x^{-1}}\dotex
x$ is the ``angular velocity in space"). We will systematically
write the left-invariant dynamics (\ref{dyng})
\begin{align}\frac{d}{dt}x(t)=DL_{ x } \omega_s(t) \end{align}

Let us suppose that  $h:G\rightarrow Y$ is a \emph{right}
equivariant smooth output map. The group action on itself by right
multiplication corresponds to the transformations $(\rho_g)_{g\in
G}$ on the output space $\YY$: for all  $x,g\in G$, $h(x
g)=\rho_g(h(x))$ i.e
$$h(R_g(x))=\rho_g(h(x))$$

Left multiplication corresponds then for the generalized body to a
change of space-fixed frame, and right multiplication to a change of
body-fixed frame. If all the measurements correspond to a part of
the state $x$ expressed in the body-fixed frame, they are affected
by a change of body-fixed frame, and the output map is right
equivariant. Thus the theory allows to build non-linear observers
such that the error equation is \emph{autonomous}, in particular for
\emph{cart-like vehicles} and \emph{rigid bodies in space}
(according to the Eulerian motion) with \emph{measurements in the
body-fixed frame} (see the example below).

\subsubsection{Observability}\label{errorsyst}

If the dimension of the output space is strictly smaller than the
dimension of the state space ($\dim y<\dim g$) the system is
necessarily not observable. This comes from the fact that, in this
case, there exists two distinct elements $x_1$ and $x_2$ of $G$ such
that $h(x_1)=h(x_2)$. If $x(t)$ is a trajectory of the system, we
have
 $\frac{d}{dt}x(t)=DL_{ g } \omega_s(t)$ and because of the
left-invariance, $g_1 x(t)$ and $g_2 x(t)$ are also trajectories of
the system:
$$
 \frac{d}{dt}(g_1  x(t))=DL_{g_1 x } \omega_s(t), \quad
 \frac{d}{dt}(g_2  x(t))=DL_{g_2 x} \omega_s(t).
$$
But since $h$ is right equivariant: $ h(g_1
x(t))=\rho_{x(t)}h(g_1)=\rho_{x(t)}h(g_2)=h(g_2 x(t)) $. The
trajectories $g_1  x(t)$ and $g_2   x(t)$ are distinct and for all
$t$ they correspond to the same output. The system is unobservable.

\subsubsection{Applying the general theory of section \ref{lie:sec}}
\label{link:sec} There are two ways to apply the theory of section
\ref{lie:sec}. i) The most natural (respecting left-invariance) does
not yield the most interesting properties: let $\UU=\RR\times\YY$
and let us look at $(u_1,u_2)=(t,h(e))$ as inputs. For all $g\in G$
let $\psi_g(t,h(e))=\left(t,\rho_{g^{-1}}(h(e))\right)$. Define a
new output map $H(x,u)=h(x)=\rho_x(h(e))=\rho_x(u_2)$. It is
unchanged by the transformation introduced in definition
\ref{inv:def} since $H(X,U)=\rho_{gx}\bigl(
\rho_{g^{-1}}(u_2)\bigr)=H(x,u)$ for all $g\in G$.
\eqref{dyng}-\eqref{output1:eq} is then a left-invariant system in
the sense of definition \ref{inv:def}, when the output map is
$H(x,u)$. ii) Let us rather look at $\omega_s(t)$ as an input :
$u(t)=\omega_s(t)\in\UU$, where $\UU=\mathfrak g\equiv\RR^n$ is the
input space. Let us define for all $g$ the map $\psi_g:
G\rightarrow\UU$ the following way
$$
\psi_g=DL_{g^{-1}}DR_g
$$
It means $\psi_g$ is the differential of the interior automorphism
of $G$. And the dynamics \eqref{dyng} writes
$\dotex x=F(x,u)=DL_x u$ and can be viewed as a right-invariant dynamics.
 For all $x,g$ we have indeed:
\begin{align*}
\dotex R_g(x)=DR_g DL_x\omega_s(t)=DL_x
DL_gDL_{g^{-1}}DR_g\omega_s(t)=DL_{R_g(x)}\psi_g(\omega_s(t))=F(R_g(x),\psi_g(u))
\end{align*}
$(\psi_g)_{g\in G}$ and $(\rho_g)_{g\in G}$ are \emph{right} group
actions
 since for all $g_1,g_2\in G$ we have
$\psi_{g_1}\circ\psi_{g_2}=\psi_{g_2g_1}$ and
$\rho_{g_1}\circ\rho_{g_2}=\rho_{g_2g_1}$.
 Thus we strictly apply the general theory of \ref{lie:sec}, exchanging the roles
of left and right multiplication.

\subsubsection{Construction of the observers}

Take $n$ linearly independent vectors $(W_1, \ldots, W_n)$ in
$TG|_e=\GG$. Consider the class of observers of the form
 \begin{align} \label{inv:obs:eq}
  \dotex{\hat x}
 &= DL_{\hat x } \omega_s(t)
  + DR_{\hat x }(\sum_{i=1}^{n}\mathcal{L}_i(\rho_{{\hat x}^{-1}}(y)) W_i)
\end{align}where the $\LL_i$'s are smooth scalar functions such that $\LL_i(h(e))=0$.
They are invariant under the transformations defined above in
section \ref{link:sec}-ii).

\subsubsection{State-error dynamics}\label{errorsystem}

The error (invariant by right multiplication) is $G\ni
\eta=(\hat{x}x^{-1})=L_{\hat{x}}({x^{-1}})$.  The error equation is
 an \emph{autonomous} differential equation \eqref{er} independent from the
trajectory $t\mapsto x(t)$ (as in the linear stationary case):
\begin{align}\label{er} \dotex \eta =DR_{\eta}( \sum_{i=1}^{n}
\mathcal{L}_i(
h(\eta^{-1}))
  W_i)
\end{align}
It can be deduced from \eqref{er:lie:eq} or directly computed using
$ \dotex \eta = DL_{\hat{x}}(\dotex
{{x}^{-1}})+D_{{x}}L_{\hat{x}}(x^{-1})\dotex {\hat{x}} $
 and
\begin{itemize}
\item $
D_{{x}}L_{\hat{x}}(x^{-1})\dotex {\hat{x}}= DR_{x^{-1}}(\dotex
{\hat{x}})=DR_{x^{-1}}DL_{\hat{x}}~\omega_s(t) +DR_{x^{-1}}DR_{\hat x} \sum_{i=1}^{n}
   \mathcal{L}_i(\rho_{{\hat x}^{-1}}(y))
  W_i
  =DR_{x^{-1}}DL_{\hat{x}}~\omega_s(t) +DR_{\eta}\sum_{i=1}^{n}
  \mathcal{L}_i(\rho_{{\hat x}^{-1}}(y))
  W_i
$
\item $
DL_{\hat{x}}(\dotex{x}^{-1})=-DL_{\hat{x}}DR_{{x}^{-1}}DL_{{x}^{-1}}\dot
x
=-DL_{\hat{x}}DR_{{x}^{-1}}\omega_s=-DR_{{x}^{-1}}DL_{\hat{x}}\omega_s(t)
$
\item $\mathcal{L}_i(\rho_{{\hat x}^{-1}}(y))=
\mathcal{L}_i(\rho_{{\hat x}^{-1}}(h(x)))
=\mathcal{L}_i(
h(\eta^{-1}))$.
\end{itemize}

\subsubsection{First order approximation}
We suppose that $\eta$ is close to $e$. Let $\xi\in\mathfrak g$ such that
 $\eta=\exp(\epsilon\xi)$ with $\epsilon\in\RR$ small. We have up to second order terms in $\epsilon$
\begin{equation*}
\dotex \xi=-\sum_{i=1}^{n} \left(\dv{\LL_i}{h}(h(e)) \dv{h}{x}(e) \xi\right) W_i
\end{equation*}

Let us define a scalar product on the tangent space $\mathfrak g$ at
$e$, and let us consider the adjoint operator of $Dh(e)$ in the
sense of the metrics associated to the scalar product. The adjoint
operator is denoted by $(Dh(e))^T$ and we take
  $
 \LL(y)=K{(Dh(e)})^T(y-h(e))
 . $
The first order approximation writes
\begin{equation}\label{er:lin:eq}
\dot\xi=-K~ Dh^T ~Dh~ \xi
\end{equation}and for $K>0$,  admits as Lyapunov function $\| \xi\|^2$
which the length of $\xi$ in the sense of the scalar product.

\subsection{A class of non-linear first-order convergent observers}

Consider for \eqref{dyng}-\eqref{output1:eq} the following observers
: $
  \dotex{\hat x}  = DL_{\hat x } \omega_s(t)
  + DR_{\hat x }[\sum_{i=1}^n[\LL_i(\rho_{\hat{x}}^{-1}(h(x)))]
  W_i]
$ where the $\LL_i$'s are smooth scalar functions such that
$\LL_i(h(e))=0$. Using the first order approximation design, take
$\mathcal{L}_1,...,\mathcal{L}_n$ such that the symmetric part (in
the sense of the scalar product chosen on $TG|_e$) of the linear map
$\xi\mapsto-\sum_{i=1}^{n} \left(\frac{\partial \LL_i}{\partial
h}({h(e))}\frac{\partial h}{\partial x}(e) \xi\right)W_i$  is
negative. When it is negative definite,  we get locally
exponentially  convergent non-linear observers around any system
trajectory.

\section{Brief example: Magnetic-aided attitude estimation}
To illustrate briefly the theory we give one of the simplest
example: magnetic-aided inertial navigation as considered in
\cite{mahony-et-al-dcd05,Bonnabel-et-al-lie-ifac08}. We just give
the system equations, the application of the theory to this example
being straightforward. It is necessary in order to pilot a flying
body to have at least a good knowledge of its orientation. This
holds for manual, or semi automatic or automatic piloting.  In
low-cost or ``strap-down" navigation systems the measurements of
angular velocity $\vw$ and acceleration $\va$ by rather cheap
gyrometers and accelerometers are completed by a measure of the
earth magnetic field $\vB$. These various measurements are fused
(data fusion) according to the motion equations of the system. The
estimation of the orientation is generally performed by an extended
Kalman filter. But the use of extended Kalman filter requires much
calculus capacity because of the matrix inversions. The orientation
(attitude) can be described by an element of the group of rotations
$\mbox{SO}(3)$, which is the configuration space of a body fixed at
a point. The motion equation are
\begin{equation}\label{dyn:eq}
    \dotex R =R(\vw\times\cdot)
\end{equation}
where
\begin{itemize}
\item $R\in \mbox{SO}(3)$ is the quaternion of norm one which represents the rotation
which maps the body frame to the earth frame,
\item $\vw$(t) is the instantaneous angular velocity vector measured
by gyroscopes and $(\vw\times\cdot)$ the skew-symmetric matrix
corresponding to wedge product with $\vw$.
\end{itemize}
If the output is the earth magnetic field  $\vec B$  measured by the
magnetometers \emph{in the body-fixed frame} $y=R^{-1}\vec B$
(\cite{bonnabel-rouchon-LN05}), the output is right equivariant. The
output has dimension $2$ (the norm of $y$ is constant) and the state
space has dimension $3$. Thus the system is not observable according
to section \ref{errorsyst}. This is why we make an additional
assumption as in
\cite{mahony-et-al-dcd05,Bonnabel-et-al-lie-ifac08}. Indeed the
accelerometers measure $\va=\dotex{\vec v}+R^{-1} \vG$ where
$\dotex{\vec v}$ is the acceleration of the center of mass of the
body and $\vG$ is the gravity vector. We suppose the acceleration of
the center of mass is small with respect to $\parallel\vG\parallel$
(quasi-stationary flight).  The measured output is thus
$y=(y_G,y_B)=(R^{-1} \vG,R^{-1} \vec B)$. One can apply the theory
as described in section \ref{link:sec}-i) or \ref{link:sec}-ii).

\section{Conclusion}

In this paper we completed the theory of \cite{arxiv-07} giving a
general framework to symmetry-preserving observers when the state
space is a Lie group. The observers are intrinsically and globally
defined. By the way, we explained the nice properties of the error
equation in two examples of \cite{arxiv-07}. In particular we
derived observers which converge around any trajectory and such that
the global behavior is independent of the trajectory as well as of
the time-varying inputs for a general class of systems.

%\bibliographystyle{plain}
%\bibliography{C:/rhn}

\begin{thebibliography}{10}

\bibitem{arnold-book-3}
V.~Arnold.
\newblock {\em Mathematical Methods of Classical Mechanics}.
\newblock Mir Moscou, 1976.

\bibitem{bonnabel:cifa06}
S.~Bonnabel, Ph. Martin, and P.~Rouchon.
\newblock Groupe de lie et observateur non-lin{\'e}aire.
\newblock In {\em CIFA 2006 (Conference Internationale Francophone
  d'Automatique), Bordeaux, France.}, June 2006.

\bibitem{Bonnabel-et-al-lie-ifac08}
S.~Bonnabel, Ph. Martin, and P.~Rouchon.
\newblock Non-linear observer on lie group for left-invariant dynamics with
  right-left equivariant output.
\newblock {\em IFAC08}, 2008.

\bibitem{arxiv-07}
S.~Bonnabel, Ph. Martin, and P.~Rouchon.
\newblock Symmetry-preserving observers.
\newblock {\em http://arxiv.org/abs/math.OC/0612193, Accepted for publication
  in IEEE AC}, Dec 2006.

\bibitem{bonnabel-rouchon-LN05}
S.~Bonnabel and P.~Rouchon.
\newblock {\em Control and Observer Design for Nonlinear Finite and Infinite
  Dimensional Systems}, chapter On Invariant Observers, pages 53--66.
\newblock Number 322 in Lecture Notes in Control and Information Sciences.
  Springer, 2005.

\bibitem{fagnani-willems-siam93}
F.~Fagnani and J.~Willems.
\newblock Representations  of symmetric linear dynamical systems.
\newblock {\em SIAM J. Control and Optim.}, 31:1267--1293, 1993.

\bibitem{grizzle-marcus-ieee85}
J.W. Grizzle and S.I. Marcus.
\newblock The structure of nonlinear systems possessing symmetries.
\newblock {\em IEEE Trans. Automat. Control}, 30:248--258, 1985.

\bibitem{hamel-mahony-icra06}
T~Hamel and R.~Mahony.
\newblock Attitude estimation on so(3) based on direct inertial measurements.
\newblock In {\em International Conference on Robotics and Automation,
  ICRA2006}, 2006.

\bibitem{koon-marsden-siam97}
W.~S. Koon and J.~E. Marsden.
\newblock Optimal control for holonomic and nonholonomic mechanical systems
  with symmetry and lagrangian reduction.
\newblock {\em SIAM J. Control and Optim.}, 35:901--929, 1997.

\bibitem{mahony-et-al-IEEE}
R.~Mahony, T.~Hamel, and J-M Pflimlin.
\newblock Non-linear complementary filters on the special orthogonal group.
\newblock {\em Accepted for publication in IEEE-AC}.

\bibitem{mahony-et-al-dcd05}
R.~Mahony, T.~Hamel, and J-M Pflimlin.
\newblock Complimentary filter design on the special orthogonal group so(3).
\newblock In {\em Proceedings of the IEEE Conference on Decision and Control,
  CDC05, Seville}, 2005.

\bibitem{martin-et-al-cocv03}
Ph. Martin, P.~Rouchon, and J.~Rudolph.
\newblock Invariant tracking.
\newblock {\em ESAIM: Control, Optimisation and Calculus of Variations},
  10:1--13, 2004.

\bibitem{morin-ieee03}
P.~Morin and C.~Samson.
\newblock Practical stabilization of driftless systems on lie groups, the
  transverse function approach.
\newblock {\em IEEE Trans. Automat. Control}, 48:1493--1508, 2003.

\bibitem{olver-book99}
P.~J. Olver.
\newblock {\em Classical Invariant Theory}.
\newblock Cambridge University Press, 1999.

\bibitem{respondek-tall-scl02}
W.~Respondek and I.A. Tall.
\newblock Nonlinearizable single-input control systems do not admit stationary
  symmetries.
\newblock {\em Systems and Control Letters}, 46:1--16, 2002.

\end{thebibliography}

                             % in the appendices.
  \end{document}